\documentclass[12pt, a4paper, reqno, oneside]{amsart}
\usepackage{amssymb}
\usepackage{amsmath}
\usepackage{amsthm}
\usepackage{graphicx}
\usepackage{color}
\usepackage[cp1251]{inputenc}
\usepackage[english]{babel}
\usepackage{cite, verbatim}

\newtheorem{theorem}{Theorem}[section]
\newtheorem{lemma}[theorem]{Lemma}
\newtheorem{proposition}[theorem]{Proposition}
\newtheorem{corollary}[theorem]{Corollary}

\newtheorem{problem}[theorem]{Problem}
\theoremstyle{definition}

\theoremstyle{remark}

\newcommand{\Aut}{\operatorname{Aut}}
\newcommand{\Out}{\operatorname{Out}}
\newcommand{\Soc}{\operatorname{Soc}}

\def\ov{\overline}
\def\wt{\widetilde}

\title[On recognition of the direct squares of simple groups]{On recognition of the direct squares of the simple groups with abelian Sylow 2-subgroups}

\author[Tao Li]{T. Li}
\address{School of Mathematics and Statistics\\  Hainan University\\
 Haikou, Hainan, P. R. China.}
\email{tli@hainanu.edu.cn}

\author[A.~R.~Moghaddamfar]{A.~R.~Moghaddamfar}
\address{Faculty of Mathematics\\ K. N. Toosi
University of Technology\\
 P. O. Box $16765$--$3381$\\ Tehran, Iran.}
\email{moghadam@kntu.ac.ir}

\author[A.~V.~Vasil'ev]{A.~V.~Vasil'ev}
\address{Sobolev Institute of Mathematics, 630090 Novosibirsk, Russia}
\email{vasand@math.nsc.ru}

\author[Zhigang Wang]{Zh. Wang}
\address{School of Mathematics and Statistics\\  Hainan University\\
 Haikou, Hainan, P. R. China.}
\email{wzhigang@hainanu.edu.cn}

\begin{document}

\thanks{A.~V.~Vasil'ev was supported by RAS Fundamental Research Program, project  FWNF-2022-0002, and by National Natural Science Foundation of China (No.~12171126).}

\begin{abstract}  The spectrum of a group is the set of orders of its elements. Finite groups with the same spectra as the direct squares of the finite simple groups with abelian Sylow 2-subgroups are considered.   It is proved that the direct square $J_1\times J_1$ of the sporadic Janko group $J_1$ and the direct squares  ${^2}G_2(q)\times{^2}G_2(q)$ of the simple small Ree groups ${^2}G_2(q)$ are uniquely characterized by their spectra in the class of finite groups, while for the direct square $PSL_2(q)\times PSL_2(q)$ of a $2$-dimensional simple linear group $PSL_2(q)$, there are always infinitely many groups (even solvable groups) with the same spectra.
 \end{abstract}

\footnotetext{{\em $2020$ Mathematics Subject Classification}: 20D60, 20D06.\\
{\em Keywords}: simple group, group with abelian Sylow 2-subgroups, small Ree group, sporadic Janko group, spectrum of a group, recognition by spectrum.}

\maketitle

\begin{center}\textsl{\small Dedicated to Professor Shi Wujie on the occasion of his 80th birthday}\end{center}

\section{Introduction}

Given a finite group $G$, the {\em spectrum} $\omega(G)$ is the set of orders of the elements of~$G$. Thirty years ago in the {\em Ricerche di Matematica}, see~\cite{BW}, R. Brandl and W. Shi established that every finite nonabelian simple group with abelian Sylow $2$-subgroups is uniquely (up to isomorphism) characterized by its spectrum in the class of finite groups. It was one of the first remarkable results in the field that has become quite popular in recent decades. In particular, the marvelous conjecture made by W. Shi in 1987 turned out to be true: every finite simple group is uniquely characterized by its spectrum and order in the class of all finite groups~\cite{VasGrMaz1}. The later result found various applications in algebra  and computational complexity, see, e.g., \cite{08KimLuRC, BP, GL}.

Here, we are interested in finite groups that are uniquely determined by their spectra only. Groups $G$ and $H$ are called {\em isospectral} if  $\omega(G)=\omega(H)$. A group $G$ is said to be {\em recognizable (by spectrum)} if $G$ is isomorphic to every group isospectral to it, {\em almost recognizable} if there are only finitely many groups isospectral to~$G$, and {\em unrecognizable} otherwise.

By the well-known theorem of V.~D. Mazurov and W. Shi \cite{MS}, a finite group which includes a nontrivial normal abelian subgroup is unrecognizable, so the problem of recognizability by spectrum is interesting only for groups with trivial solvable radical. Despite the fact that there are a huge number of simple groups recognizable by spectrum (see \cite[Theorem 2.1]{GMSVY}), up to now there were only two examples of a recognizable group that is the direct square of a simple group: $Sz(2^7)\times Sz(2^7)$ \cite{M} and $J_4\times J_4$ \cite{GM}. The reason for this is that the standard methods of proving the nonsolvability of groups isospectral to a simple group are based on the properties of its prime graph. Recall that the prime graph $\Gamma(G)$ of a group $G$ is a graph with vertex set $\pi(G)$, the set of prime divisors of the order of $G$, in which two vertices $p$ and $q$ are adjacent if and only if $p\neq q$ and $pq\in\omega(G)$; we also denote by $t(G)$ the maximum number of pairwise nonadjacent vertices in~$\Gamma(G)$. Clearly, such methods are not applicable to proving the nonsolvability of a group isospectral to the square of a group, since in this case, the corresponding graph is always complete.

Recently in \cite[Theorem~1]{WVGZ}, a helpful criterion of nonsolvability of a finite group was established. A corollary of this criterion yields that for every finite group $L$ with $t(L)\geq 4$, a group $G$ isospectral to the direct square $L\times L$ must be nonsolvable \cite[Corollary 2]{WVGZ}. Using this result, the authors of that paper proved that the direct square $Sz(q)\times Sz(q)$ of the simple exceptional Suzuki group $Sz(q)$ is recognizable by spectrum for every $q=2^\alpha$, $5\neq\alpha\geq 3$ odd, while for $Sz(2^5)\times Sz(2^5)$, there are exactly four finite groups with the same spectra \cite[Theorems 3 and 4]{WVGZ}.

In the present paper, we consider the direct squares of the finite nonabelian simple groups with abelian Sylow $2$-subgroups. The finite groups with abelian Sylow $2$-subgroups were described (up to some details covered only by the classification of finite simple groups) by J.~H. Walter in \cite{Wa}, see details in Lemmas~\ref{lem1} and~\ref{l:Walter} below. Essentially, a nonabelian simple group with abelian Sylow $2$-subgroups is either a $2$-dimensional linear group $L_2(q)=PSL_2(q)$ (here and further we use the {\em Atlas one-letter} notation for simple groups, see \cite{Atlas}), a small Ree group $R(q)={^2}G_2(q)$, $q=3^\alpha$, $\alpha\geq 3$ odd, or the sporadic Janko group $J_1$.

It is well known that $t(L)=3$ for $L=L_2(q)$, $t(L)=4$ for the Janko group $J_1$, and $t(L)=5$ for the small Ree groups $R(q)$~\cite{VasVd05}. It turns out that the direct square of a simple group $L=L_2(q)$ (regardless of whether its Sylow $2$-subgroups are abelian or not) is unrecognizable by spectrum, moreover, it is easy to construct a solvable group isospectral to $L$, see Proposition~\ref{pro1} below. On the other hand, all other nonabelian simple groups with abelian Sylow $2$-subgroups are recognizable by spectrum.

\begin{theorem}\label{th1} Let $L=R(q)$, where $q=3^{\alpha}$, $\alpha\geq 3$ odd. If $G$ is a finite group isospectral to $L\times L$, then $G\simeq L\times L$.
\end{theorem}

\begin{theorem}\label{th2} If $G$ is a finite group isospectral to $J_1\times J_1$, then $G\simeq J_1\times J_1$.
\end{theorem}

As a consequence of these results, we completely solve the recognition problem for the direct squares of the simple groups with abelian Sylow $2$-subgroups.

\begin{corollary}\label{cor} Suppose that $L$ is a simple group with abelian Sylow $2$-subgroups. Then the direct square $L\times L$ of $L$ is recognizable by spectrum if $L$ is either a small Ree group $R(q)$ or the sporadic Janko group $J_1$, and $L\times L$ is unrecognizable otherwise.
\end{corollary}

Walter's classification and the results obtained here lead naturally to the next challenging question.

\begin{problem} Which finite groups with abelian Sylow $2$-subgroups are recognizable by their spectra?
\end{problem}

The rest of the paper is organized as follows. The preliminary general results are collected in Section~2. The necessary information on the groups with abelian Sylow $2$-subgroups is listed in Section~3. Sections 4 and 5 contain the proofs of Theorems~\ref{th1} and~\ref{th2} respectively.

\section{Preliminaries}

For a natural number~$n$, we write $\pi(n)$ to denote the set of prime divisors of~$n$. If $\pi$ is the set of primes, then $\pi'$ is the set of primes not in~$\pi$.

Let $G$ be a group. Clearly, the spectrum $\omega(G)$ of $G$ includes the set $\pi(G)=\pi(|G|)$ of the prime divisors of~$|G|$, and the prime graph $\Gamma(G)$ is determined by~$\omega(G)$. Furthermore, $\omega(G)$ is closed with respect to taking divisors and therefore it is uniquely defined by its subset $\mu(G)$ consisting of elements maximal with respect to the divisibility relation.

As usual, $\Aut G$ and $\Out G$ are respectively the automorphism group and the outer automorphism group of~$G$, and $\Soc(G)$ is its socle, i.e., the product of the minimal normal subgroups. If a group $H$ acts on a group $K$, then $K\rtimes H$ denotes their natural semidirect product.

\begin{lemma}[{Bang \cite{86Bang}, Zsigmondy \cite{Zs}}]\label{l:bangz} Let $q, n\geq 2$ be integers. Then either there is a prime number $r$ that divides $q^n-1$ and does not divide $q^i-1$ for all $i<n$, or one of the following conditions is satisfied:
\begin{enumerate}
 \item[(1)] $q=2$ and $n=6$;
 \item[(2)] $q$ is a Mersenne prime and $n=2$.
 \end{enumerate}
\end{lemma}

A prime number $r$ from Lemma~\ref{l:bangz} is called a \emph{primitive prime divisor} of $q^n-1$. Observe that $n$ is the multiplicative order of $q$ modulo $r$.

Primitive prime divisors are a valuable tool when dealing with element orders of groups of Lie types. If $G$ is such a group over the field of order $q$ and is not a Suzuki or Ree group, then the adjacency of a prime $r$ in $\Gamma(G)$ depends mostly not on $r$ itself but on the order of $q$ modulo $r$, see details in \cite{VasVd05,VV}. Unfortunately, the situation for Suzuki and Ree groups is a bit more complicated. We overcome this in the case of small Ree groups with the help of the following lemma, which is a particular case of the main result of \cite{G}.

\begin{lemma}\label{l:GP}
For every odd integer $n\geq 3$ and each $\varepsilon\in\{+, -\}$, there is a prime $r$ that divides $3^n+\varepsilon3^{(n+1)/2}+1$ and does not divide $3^{2i}-3^{i}+1$ for all odd $i<n$.
\end{lemma}

Analogously, we refer to a number $r$ from Lemma~\ref{l:GP} as a \emph{primitive prime divisor} of $3^n+\varepsilon3^{(n+1)/2}+1$. Note that $3^{2n}-3^n+1=(3^n-3^{(n+1)/2}+1)(3^n+3^{(n+1)/2}+1)$ and the union of all primitive prime divisors of $3^n+\varepsilon 3^{(n+1)/2}+1$, $\varepsilon\in\{+, -\}$, is the set of primitive prime divisors of $3^{6n}-1$.

The next lemma and its corollary are the essential ingredients of our proofs.

\begin{lemma}\label{l:criterion}
Let $G$ be a finite group. Suppose that there is a subset $\sigma(G)$ of $\pi(G)$
such that one of the following holds:
\begin{enumerate}
\item[(1)] $|\sigma(G)|=3$ and $pq\not\in\omega(G)$ for all distinct $p, q\in\sigma(G);$
\item[(2)] $|\sigma(G)|=4$, $pq\in\omega(G)$ for all distinct $p, q\in\sigma(G)$, but $pqr\not\in\omega(G)$ for all pairwise distinct $p, q, r\in \sigma(G).$
\end{enumerate}
Then $G$ is nonsolvable.
\end{lemma}

\begin{proof} If (1) holds, then it follows from \cite[Theorem~1]{57Hig}. For (2), it is \cite[Theorem~1]{WVGZ}.
\end{proof}

\begin{corollary}\label{cor1}  Let $L$ and $G$ be finite groups. Suppose that one of the following holds:
\begin{enumerate}
\item[(1)] $t(L)\geq 3$ and $\Gamma(G)=\Gamma(L);$
\item[(2)] $t(L)\geq 4$ and $\omega(G)=\omega(L\times L)$.
\end{enumerate}
Then $G$ is nonsolvable.
\end{corollary}

We need the following two easily established facts.

\begin{lemma}\label{l:cyclic1} If $p$ and $q$ are odd prime divisors of the order of a solvable group $G$, nonadjacent in $\Gamma(G)$, then either a Sylow $p$-subgroup or a Sylow $q$-subgroup of $G$ is cyclic.
\end{lemma}

\begin{proof} We may assume that the group $G$ includes a $\{p, q\}$-subgroup which is a extension of a $p$-group by a Sylow $q$-group of~$G$. The rest follows from \cite[Theorem~3.16]{Gor}.
\end{proof}

\begin{lemma}\label{l:cyclic2} Suppose that for a prime~$p$, a Sylow $p$-subgroup of the solvable radical $R(G)$ of a nonsolvable group $G$ has a unique subgroup $P$ of order $p$. If $C=C_G(P)$ is the centralizer of $P$ in $G$ and $R(C)$ is the solvable radical of $C$, then $C$ has the same nonabelian factors as $G$ and $\Soc(C/R(C))\simeq\Soc(G/R(G))$.
\end{lemma}

\begin{proof} On the one hand, the normalizer $N=N_G(P)$ of $P$ in $G$ includes the normalizer of a Sylow $p$-subgroup of $K$ containing~$P$. Therefore, by the Frattini argument, $G=NR(G)$, so $N$ has the same nonabelian factors as $G$ and $\Soc(N/R(N))\simeq\Soc(G/R(G))$, where $R(N)$ is the solvable radical of~$N$. On the other hand, the quotient $N/C$, being isomorphic to an automorphism group of the cyclic group~$P$, must be abelian. Thus, $C$ has the same nonabelian factors as $G$ and $\Soc(C/R(C))\simeq\Soc(G/R(G))$, as required.
\end{proof}

Let $\pi$ be a set of primes. Following~\cite{56Hal},  we say that a finite group $G$ is a $D_\pi${\em-group} (or $G$ has the property $D_\pi$) if $G$ contains a Hall $\pi$-subgroup, all Hall $\pi$-subgroups of $G$ are conjugate, and every $\pi$-subgroup of $G$ is contained in some Hall $\pi$-subgroup of $G$.

\begin{lemma}{{\rm\cite[Theorem~6.6]{11VR}}}\label{l:dpi}
Let $\pi$ be a set of primes, let $G$ be a group and $N$ a normal subgroup of $G$. The group $G$ is a $D_\pi$-group if and only if $N$ and $G/N$ are $D_\pi$-groups.
\end{lemma}

A group $G$ is called a {\em cover} of a group $H$ if there exists an epimorphism from $G$ onto~$H$. The following two lemmas contain the well-known facts about the orders of elements in a cover of a finite group.

\begin{lemma}\label{l:frob}
Given a prime $p$, suppose that $V$ is a normal abelian $p$-subgroup of a finite group $G$ and $G/V$ is a Frobenius group with a kernel $N$ and a cyclic complement~$\langle x\rangle$. If $p$ does not divide $|N|$ and $N\not\subseteq C_G(V)/V$, then $C_V(x)\neq 1$ and $G$ contains an element of order~$p|x|$.
\end{lemma}

\begin{proof} See, e.g., \cite[Lemma~1]{94Maz.t}.
\end{proof}

\begin{lemma}\label{l:hh} Given an odd prime $p$, suppose that $V$ is a normal $p$-subgroup of a finite group $G$ and $G/V$ is a semidirect product of a $p'$-group $N$ with abelian Sylow $2$-subgroups and a cyclic $p$-group~$\langle x\rangle$. If $O_p(G/V)=1$ and $N\not\subseteq VC_G(V)/V$, then $G$ contains an element of order~$p|x|$.
\end{lemma}

\begin{proof}
It follows from the Hall-Higman theorem~\cite[Theorem~B]{56HalHig}.
\end{proof}

\section{Simple groups with abelian Sylow 2-subgroups}

Walter~\cite{Wa} obtained a characterization of finite groups with abelian Sylow $2$-subgroups (up to some details covered only by the classification of finite simple groups).

\begin{lemma}\label{lem1}  If $G$ is a nonabelian simple group with an abelian Sylow $2$-subgroup $S$, then $S$ is elementary abelian and $G$ is one of the following groups:
\begin{itemize}
\item[{(1)}] a linear group $L_2(q)$, $q>3$, $q\equiv 3, 5 \pmod{8}$ or $q=2^\alpha;$
\item[{(2)}] a small Ree group $R(q)$, $q=3^\alpha$,  $\alpha\geq 3$ is odd{\rm;}
\item[{(3)}] the sporadic Janko group $J_1$.
\end{itemize}
\end{lemma}

\begin{lemma}\label{l:Walter}{\em\cite[Theorem~I]{Wa}} If $G$ is a finite group with an abelian Sylow $2$-subgroup, then the group $O^{2'}(G/O_{2'}(G))$ is a direct product of an abelian $2$-group and some nonabelian simple groups listed in Lemma~{\em\ref{lem1}}.
\end{lemma}

The following lemma collects some basic information on the simple linear groups~$L_2(q)$.

\begin{lemma} \label{l:lin} Let $G=L_2(q)$, where $q=p^\alpha\geq 4$, and $d=(2, q-1)$. Then the following hold:
\begin{enumerate}
 \item[(1)] $|G|=q(q^2-1)/d;$
 \item[(2)] $\mu(G)=\{p, (q-1)/d, (q+1)/d\};$
 \item[(3)] if $G\leq H\leq \Aut G$ and $|H:G|$ is odd, then  $H=G\rtimes\langle\psi\rangle$, where $\psi$ is a field automorphism of $G$ and $C_G(\psi)\simeq L_2(q^{1/|\psi|});$
 \item[(4)] if $q>9$, then the Schur multiplier of $G$ is of order~$d$;
 \item[(5)] $G$ includes a Frobenius subgroup with kernel of order $q$ and complement of order $(q-1)/d$.
\end{enumerate}
\end{lemma}

We are ready to prove that there are infinitely many finite groups isospectral to the direct square of the group $L_2(q)$ for every $q\geq 4$.

\begin{proposition}\label{pro1} If $L=L_2(q)$,  where $q=p^\alpha\geq 4$, $p$ a prime, then the direct square $L\times L$ is unrecognizable by spectrum.
\end{proposition}
\begin{proof}
Let $d=(2,q-1)$. We write $E_1$ and $E_2$ for the elementary abelian $p$-groups of orders $q$ and $q^2$; and $C_1$ and $C_2$ for the cyclic subgroups of orders $(q-1)/d$ and $(q+1)/d$, respectively. Let $F_1$ and $F_2$ be Frobenius groups with kernels $E_1$ and $E_2$, and complements $C_1$ and $C_2$, respectively. Then $\mu(F_1\times F_2)=\{p(q-1)/d, p(q+1)/d, (q^2-1)/d^2\}=\mu(L\times L)$ in view of Lemma~\ref{l:lin}(2). Now $L\times L$ is unrecognizable by the main result of~\cite{MS}.
\end{proof}

The following two lemmas help to deal with the spectra of covers of groups~$L_2(q)$.

\begin{lemma}\label{l:LinRep} Let $G=L_2(q)$, $q=p^\alpha\geq 4$, $p$ a prime and $\alpha>1$, act on $\mathbb{F}_r$-module $V$, where $r\neq p$ is an odd prime. Then $C_V(x)>0$ for every $p$-element $x$ of $G$.
\end{lemma}

\begin{proof} If the action of $G$ on $V$ is not faithful, then $C_G(V)=G$, and we are done. If $G$ acts faithfully, then the lemma follows, e.g., from \cite[Proposition~1.2]{08DiMZal}.
\end{proof}

\begin{lemma}\label{l:LinAct}
Let $K$ be a normal subgroup of a finite group $G$ such that $G/K\simeq L_2(q)$, where $q=p^\alpha\geq 4$, $q\neq 9$, $p$ a prime, and $(|K|,2p)=1$. If $x$ and $y$ are nontrivial $p$-elements of $G$, then $\pi(C_G(x))=\pi(C_G(y))$.
\end{lemma}

\begin{proof} Since $p$ does not divide $|K|$, we have $|x|=|y|=p$. By the Sylow theorem, we may also assume that $x$ and $y$ lie in the same Sylow $p$-subgroup $P$ and their images $\ov{x}$ and $\ov{y}$ in $\ov{G}=G/K$ lie in the image $\ov{P}$ of $P$ isomorphic to~$P$. If $P$ (and so $\ov{P}$) is cyclic, then there is nothing to prove, so we may assume that $q>p$. In this case, we will prove that $\pi(C_G(y))=\pi(K)$ for every $p$-element $y$.

Suppose to the contrary that $r\in\pi(K)\setminus\pi(C_G(y))$. Since $|K|$ is odd, $K$ is solvable and, in particular, $r$-solvable. Consider the {\em upper $r$-series} of $K$ that is a normal $r$-series:
\begin{equation}\label{eq:r}
1=R_0\leq K_1<R_1\leq K_2\leq\ldots\leq R_{t-1}<K_t\leq R_t=K,
\end{equation}
where $K_i/R_{i-1}=O_{r'}(K/R_{i-1})$ and $R_i/K_i=O_r(K/K_i)$ for $i=1, \ldots, t$.

First we suppose that $R=R_t/K_t\neq 1$. Put $V=R/\Phi(R)$ for the factor group of $R$ by its Frattini subgroup $\Phi(R)$. Since $\ov{G}$ acts on $V$, we arrive at a contradiction by Lemma~\ref{l:LinRep}. Thus, $K=K_t>R_{t-1}$ and $R_{t-1}>K_{t-1}$. It suffices to prove that $\widetilde{K}=K/R_{t-1}$ must be a direct factor in $\widetilde{G}=G/R_{t-1}$, because if so, $\ov{G}$ acts on $V=R/\Phi(R)$, where $R=R_{t-1}/K_{t-1}$ and we can argue as in the previous case.

We will prove by induction on the order of $\widetilde{K}$ that $\widetilde{G}=\widetilde{K}\times\widetilde{L}$, where $\widetilde{L}\simeq L_2(q)$. Let $U$ be a minimal normal in $\wt{G}$ subgroup of $\wt{K}$. Then $U$ is an elementary abelian $s$-subgroup for some odd prime $s\neq p$. By inductive hypothesis, $\wt{G}/U=\wt{K}/U\times\wt{H}/U$, where $\widetilde{H}/U\simeq L_2(q)$. If $C_{\wt{H}}(U)=U$, then $[U,\wt{y}]\neq 1$, where $\wt{y}$ is the image of $y$ in $\wt{G}$. It follows that $F=[U,\wt{y}]\langle\wt{y}\rangle$ is a Frobenius group acting on the $r$-group~$V$. Moreover, $C_{K/K_{t-1}}(V)\leq V$ due to definition of the upper $r$-series~\eqref{eq:r}. The application of Lemma~\ref{l:frob} yields $C_V(\wt{y})\neq 1$, so $r\in\pi(C_G(y))$, a contradiction.

Thus, $C_{\wt{H}}(U)=\wt{H}$. Since the Schur multiplier of $L_2(q)$, $q\neq 9$, is a $2$-group, it follows that $\wt{H}=U\times\wt{L}$, where $\wt{L}\simeq L_2(q)$. Since $[\wt{K},\wt{L}]\leq U$, the equalities $[\widetilde{K},\widetilde{L},\widetilde{L}]=1$ and $[\widetilde{L},\widetilde{K},\widetilde{L}]=1$ hold, and hence $[\widetilde{L},\widetilde{K}]=[\widetilde{L},\widetilde{L},\widetilde{K}]=1$.
\end{proof}

The next two lemmas accumulate information on the small Ree groups $R(q)={}^2G_2(q)$, that we need in further considerations.

\begin{lemma} \label{l:ree} Let $G=R(q)$, where $q=3^\alpha$, $\alpha\geq 3$ is odd. Then the following hold:
\begin{enumerate}
 \item[(1)] $|G|=q^3(q-1)(q^3+1);$
 \item[(2)] $\mu(G)=\{6, 9, q-1, (q+1)/2, q-\sqrt{3q}+1, q+\sqrt{3q}+1\};$
 \item[(3)] $\Aut G=G\rtimes \langle \varphi\rangle$, where $\varphi$ is a field automorphism of $G$ of order~$\alpha$, and if $\psi\in\langle\varphi\rangle$, then $C_G(\psi)\simeq R(q^{1/|\psi|});$
 \item[(4)] the Schur multiplier of $G$ is trivial.
\end{enumerate}
\end{lemma}

\begin{proof} Items (1) and (3) were proved by Ree~\cite{R}. Item (2) is \cite[Lemma~4]{BW}. Item (4) was established in~\cite{66AlpGor}.
\end{proof}

\begin{lemma}\label{l:ReeRep}
Let $G=R(q)$, where $q\geq 27$, and  $g\in G$ an $r$-element for a prime~$r\neq 3$. If $G$ acts faithfully on a $p$-group $V$ for some prime $p$, then the coset $Vg$ of the natural semidirect product $V\rtimes G$ contains an element of order~$p|g|$.
\end{lemma}

\begin{proof}
We can suppose that $V$ is an elementary abelian $p$-group and consider $V$ as an irreducible $G$-module. If $p\neq 3$, then in virtue of \cite[Theorem~1.1]{TiepZal}, the minimal polynomial of $g$ in this representation is equal to $x^{|g|}-1$. Hence, there exists $v\in V$ such that $v(1+g+g^2+\dots+g^{|g|-1})\neq 0$. Then $(vg)^{|g|}\neq 0$, so $|vg|=p|g|$. On the other hand, by \cite{GurT}, the group $G$ is unisingular, that is every element of $G$ has a nonzero fixed point in every $G$-module over a field of characteristic~$3$, so the lemma holds in the case $p=3$.
\end{proof}

Now we turn to the sporadic Janko group $J_1$.

\begin{lemma}\label{spectraJ1} Let $G=J_1$. Then the following hold:
\begin{enumerate}
\item [(1)] $|G|=2^3\cdot 3\cdot 5\cdot 7\cdot 11\cdot 19;$
\item [(2)] $\mu(G)=\{6, \ 7, \ 10, \  11, \ 15, \  19\};$
\item [(3)] $\Aut G=G$ and the Schur multiplier of $G$ is trivial.
\end{enumerate}
\end{lemma}

\begin{proof} See \cite[p.~36]{Atlas}
\end{proof}

\begin{lemma}\label{l:J1Rep}
Let $G=J_1$ and let $g\in G$ be an element of odd prime order $p\neq 11$. If $G$ acts faithfully on a $p$-group $V$, then the coset $Vg$ of the natural semidirect product $V\rtimes G$ contains an element of order~$p^2$.
\end{lemma}

\begin{proof}
We can suppose that $V$ is an elementary abelian $p$-group and consider $V$ as an absolutely irreducible $G$-module. The rest follows from~\cite[Theorem~1.1]{99Zal}.
\end{proof}

The following lemmas list some general properties of simple groups with abelian Sylow $2$-subgroups.

\begin{lemma}\label{NorSyl3} If $G$ is a nonabelian simple group with an abelian Sylow $2$-subgroup, then the normalizer of a Sylow $3$-subgroup of $G$ is solvable.
\end{lemma}

\begin{proof} Let $T$ be a Sylow $3$-subgroup of $G$ and $N$ the normalizer of $T$ in $G$. According to Lemma~\ref{lem1}, $G$ is isomorphic to $L_2(q)$, $R(q)$, or $J_1$. If $G=J_1$, then $N$ is isomorphic to the direct product of two dihedral groups of orders $6$ and $10$ \cite[p.~36]{Atlas}. Hence we may assume that $G$ is a group of Lie type. If the characteristic of $G$ equals $3$, then $N$ is a Borel subgroup of $G$, and hence $N$ is solvable. Thus, it remains to deal with the case, when $G=L_2(q)$ and $3$ divides the order of some torus $U$ of $G$. Here, on the one hand, $N$ includes the normalizer $M$ of $U$ in $G$, while on the other hand, $M$ is a maximal subgroup of $G$, so $N=M$ is isomorphic to a dihedral group (see, e.g., \cite[Table~8.1]{Bray}).
\end{proof}

\begin{lemma}\label{lem2} Let $G$ be a nonabelian simple group with an abelian Sylow $2$-subgroup. Then $|G|$ is coprime to $5$ if and only if $G$ is one of the following groups:
\begin{itemize}
\item[{(1)}] $L_2(q)$ and $q\equiv 2, 3 \pmod{5};$
\item[{(2)}] $R(q)$, $q=3^{\alpha}$, $\alpha\geq 3$ is odd.
\end{itemize}
\end{lemma}
\begin{proof} This can be easily deduced from Lemmas~\ref{l:lin}(2) and~\ref{l:ree}(2) and the fact that $5$ divides $p^4-1$ whenever $p$ is a prime distinct from~$5$.
\end{proof}

\section{Proof of Theorem \ref{th1}}

Let $L=R(q)$, where $q=3^\alpha$, $\alpha\geq 3$ odd. By Lemma \ref{l:ree}(2),
$$\mu(L)=\{6, \ 9, \ q-1, \  (q+1)/2, \ q-\sqrt{3q}+1, \  q+\sqrt{3q}+1\}.$$
 Put $m_1=6$, $m_2=9$, $m_3=q-1$, $m_4=(q+1)/2$,  $m_5=q-\sqrt{3q}+1$, $m_6=q+\sqrt{3q}+1$. As easily seen, if $1\leq i< j\leq 6$, then $(m_i, m_j)=2$ for $i, j\in\{1,3,4\}$ and $(m_i,m_j)=1$ otherwise. In particular, $t(L)=5$.  Suppose $G$ is a finite group such that $$\mu(G)=\mu(L\times L)=\left\{[m_i,  m_j]  \  |  \ 1\leq i<j\leq 6,\ (i, j)\neq (1,3), (1,4)\right\}.$$
Corollary \ref{cor1} implies that $G$ is nonsolvable.

Let $K$ be the solvable radical of $G$ and $\overline{G}=G/K$. By Lemma~\ref{l:Walter}, the socle $S=\Soc(\ov G)$ of $\overline{G}$, i.e., the product of the minimal normal subgroups of $G$, is a direct product of finite nonabelian simple groups: $S=L_1\times L_2\times \cdots \times L_k$, where factors $L_i$ are isomorphic to some groups from the conclusion of Lemma~\ref{lem1}. Since $5\not\in\omega(G)$, Lemma \ref{lem2} makes the list of possible simple factors $L_i$ even narrower.

\begin{lemma}\label{PosLi}
For every $i=1,\ldots,k$, either $L_i\simeq R(u)$, $u=3^{\alpha_i}$, $\alpha_i\geq 3$ odd, or $L_i\simeq L_2(u)$, where $u\equiv 2, 3\pmod{5}$ and either $u\equiv 3, 5\pmod{8}$ or $u=2^{\alpha_i}$, $\alpha_i\geq 3$ odd. In particular, if $L_i\not\in\{L_2(8), L_2(13)\}$, then $u\geq 27$.
\end{lemma}

Obviously, $2,3\in\pi(L_i)$ for every $L_i, i=1, \ldots, k$.

Put $\pi_1=\{2\}$, $\pi_2=\{3\}$, and $\pi_i=\pi(m_i)\setminus\{2\}$ for $i=3,4,5,6$. Further, we will often use the obvious fact that the product of three odd primes from different $\pi_i$ cannot belong to $\omega(G)$.

\begin{lemma}\label{l:3dif}
For every three primes from different $\pi_i$, $i\in\{3,4,5,6\}$, at least two of them divide $|\ov{G}|$. In particular, there are distinct $a,b\in\{3,4,5,6\}$ such that $\pi_a\cup\pi_b\subseteq\pi(\ov{G})$.
\end{lemma}

\begin{proof} Suppose to the contrary that three primes $p, q, s$ from different $\pi_i$, $i\in\{3,4,5,6\}$, divide only $|K|$. Let $R$ be a Sylow $r$-subgroup of $G$ for a prime $r$ from the remaining fourth component. Observe that $K$ and $\ov{G}$ are $D_\sigma$-groups for $\sigma=\{p, q, s, r\}$, because $K$ is solvable and $\sigma\cap\pi(\ov{G})=\{r\}$. It follows from Lemma~\ref{l:dpi} that $G$ is also a $D_\sigma$-group. Consider now a Hall $\sigma$-subgroup $H$ of the solvable group $KR$. It is clear that $H$ is a Hall $\sigma$-subgroup of~$G$. Since every $\sigma$-subgroup of $G$ is conjugate to a subgroup of $H$ (cf. the definition of a $D_\pi$-group), the spectrum of $H$ contains the products of every two distinct primes from $\sigma$ and does not contain the product of any three of such primes. Thus, the group $H$ and the set $\sigma$ satisfy the hypothesis of Lemma~\ref{l:criterion}(2); this contradicts the solvability of~$H$.
\end{proof}

Further we denote by $G_1$ the preimage of $S$ in $G$.

\begin{lemma}\label{kLessThan3}
The number $k$ of factors of $S$ does not exceed $2$.
\end{lemma}

\begin{proof} Suppose to the contrary that $k\geq 3$. First, assume that $\pi(S)$ meets nontrivially at least two of $\pi_i$ for $i\geq 3$. Then, having in mind that the order of every $L_i$ must be divided by a prime greater than $3$, we may assume that up to reordering of direct factors in $S$, $\pi(L_1)$ contains a prime $p$ from one of these intersections and $\pi(L_2)$ contains a prime $q$ from another. Since $k\geq 3$, it follows that $3pq\in\omega(G)$; a contradiction.

Thus, there is at most one index $j$ from $\{3,4,5,6\}$ such that $\pi_j\cap\pi(S)\neq \varnothing$. Put $\sigma=\pi(G)\setminus(\pi_j\cup\{2\})$. Let $T$ be a Sylow $3$-subgroup of $S$ and $T_i=T\cap L_i$, $i=1, \ldots, k$. Lemma~\ref{NorSyl3} implies that the normalizer $N_{L_i}(T_i)$ is solvable for all $i$. In view of Lemma~\ref{l:Walter}, the quotient $G/G_1$ is also solvable as a group of odd order. Thus, the same holds true for $\ov{N}=N_{\ov{G}}(T)$ and hence for its preimage $N$ in~$G$.

Consider a Hall $\sigma$-subgroup $H$ of $N$. By the Frattini argument, $\ov{N}S=\ov{G}$, so $H$ is a Hall $\sigma$-subgroup for the whole group $G$. On the other hand, $G$ is a $D_\sigma$-subgroup, because $K$ and $G/G_1$ are solvable, while $\pi(S)\cap\sigma=\{3\}$. Arguing as in the proof of the previous lemma, we see that the solvable group $H$ satisfies the hypothesis of Lemma~\ref{l:criterion}(2), a contradiction.
\end{proof}

\begin{lemma}\label{k=1andL=L2}
If $k=1$, then $S$ is a Ree group.
\end{lemma}

\begin{proof} Suppose to the contrary that $S$ is not a Ree group. Then $S=L_1\simeq L_2(u)$, $u=v^{\alpha_1}$, and $u$ satisfies the conditions of Lemma~\ref{PosLi}. Observe that if $G/G_1\neq 1$, then $\ov{G}$ is a split extension of $S$ by a field automorphism of odd order in view of Lemmas~\ref{l:Walter} and~\ref{l:lin}(3).

Suppose that there is an odd prime $r\in\pi(G)$ such that $r$ divides only the order of~$S$. If $x\in G$ has order $r\neq v$, then $\pi(C_G(x))=\pi(G)$, because a Sylow $r$-subgroup of $G$ is cyclic, while $r$ is adjacent to any prime from $\pi(G)\setminus\{r\}$ in the prime graph~$\Gamma(G)$. If $r=v$, then in view of Lemma~\ref{l:LinAct}, the same equality $\pi(C_G(x))=\pi(G)$ holds for some element $x$ of order~$r$ (take an element $x$ whose image $\ov{x}$ in $\ov{G}$ lies in the centralizer of a field automorphism of $S$ generating $\ov{G}/S$ and apply Lemma~\ref{l:lin}(3)).  Since $\pi(C_G(x))=\pi(G)$, we may take three primes $p, q, s$ such that they together with $r$ lie in four distinct components $\pi_i, i\in\{2, 3, 4, 5, 6\}$. Since the image of $C_G(x)$ in $\ov{G}$ lies in $C_{\ov{G}}(\ov{x})$ which is solvable, it follows that $C_G(x)$ is also solvable. Then, by Lemma~\ref{l:criterion}(1), at least two of $p, q, s$ must be adjacent in~$\Gamma(C_G(x))$, let, say, $pq\in\omega(C_G(x))$. Then $rpq\in\omega(C_G(x))\subseteq\omega(G)$, which is impossible. Thus, every odd prime from $\pi(G)$ must divide $|K|\cdot|G/G_1|$.

The arguments similar to the ones from the previous paragraph show that for every odd prime $r\in\pi(S)$, a Sylow $r$-subgroup $R$ of $G$ cannot be cyclic. Indeed, in this case it is easy to see that $\pi(C_G(R))=\pi(G)$. So it suffices to prove that $C_G(R)$ is solvable. Observe that $R$, being cyclic, cannot intersect $G/G_1$ (cf., e.g., \cite[Lemma~3.10]{15Vas}). Hence $C_G(R)$ is solvable, because its image in $\ov{G}$ is solvable.

Let $\sigma=\pi(u(u-1))\setminus\{2\}$. We claim that $S$ has the property $D_\sigma$, and hence $G$ also has this property. Indeed, if $u$ is a power of $2$ or $u\equiv3\pmod8$, then a Borel subgroup is a required Hall $\sigma$-subgroup of $S$, while if $u\equiv3\pmod8$, then so is the subgroup of index $2$ in a Borel subgroup of $S$. In either case, the normalizer of this Hall $\sigma$-subgroup is a Borel subgroup. Consider the preimage $H$ of this normalizer in $G$. It is a solvable subgroup, so by Lemma~\ref{l:criterion}(2), $\pi(H)$ meets nontrivially at most three different components $\pi_i$, $i=2,\ldots,6$. On the other hand, $|G:H|=u+1$, hence the set $\pi(u+1)$ must include at least two of these components. Furthermore, $S$ and, therefore, $G$ contains an element of order $(u+1)/(2,u+1)$. Thus, $\pi(u+1)\setminus\{2\}=\pi_a\cup\pi_b$ for some distinct $a,b\in\{2,\ldots,6\}$. It follows that the set $\pi(u(u-1))\setminus\{2\}$ contains primes only from three other components and, by similar arguments involving a Hall $(\pi(u+1)\setminus\{2\})$-subgroup of $S$, includes at least two of these components.

The solvable radical $K$ of $G$ has a normal $2$-complement $C=O_{2'}(K)$ due to Lemma~\ref{l:Walter}.
Since $G/G_1$ is cyclic, $\pi(G/G_1)$ meets nontrivially at most two of the odd components $\pi_i$.
It follows that $\pi(C)$ includes at least three odd components~$\pi_i$.

Suppose first that $\pi(u-1)$ meets two odd components, say $\pi_c$ and $\pi_d$. Obviously, $\{a, b\}\cap\{c, d\}=\varnothing$. Lemma~\ref{l:cyclic2} implies now that $C$ cannot have cyclic Sylow subgroups. If $\pi(C)$ contains primes from four different odd components $\pi_i$, then we immediately derive to contradiction applying Lemmas~\ref{l:cyclic1} and~\ref{l:criterion}(2). Let $\pi(C)$ be the union of three components. Clearly, there is an odd prime $r\in\pi(G)\setminus\pi(C)$ dividing the order of $S$. By the preceding arguments, a Sylow $r$-subgroup $R$ of $G$ is not cyclic. Then the group $H=CR$ is solvable and all its Sylow subgroups are not cyclic, which is again impossible in view of Lemmas~\ref{l:cyclic1} and~\ref{l:criterion}(2).

Thus, $\pi(u-1)$ meets nontrivially only one odd component. It follows that $\pi(u-1)\setminus\{2\}=\pi_c$ and $\{v\}=\pi_d$ for some $c, d\in\{2, \ldots, 6\}\setminus\{a, b\}$. In particular, $v$ is odd. In this case, $\pi(G/G_1)$ cannot contain the primes from two different odd components (recall that $G/G_1$ is generated by a field automorphism). Therefore, $\pi(C)$ includes at least four odd components~$\pi_i$. Lemma~\ref{l:cyclic2} implies that a Sylow $p$-subgroup of $C$ is not cyclic for each $p\in\pi(C)\setminus(\pi_a\cup\pi_b)$.
If $\pi(C)$ meets nontrivially all three odd components distinct from $\pi_a$ and $\pi_b$, then we get a contradiction as at the end of the previous paragraph, adding to $C$ a Sylow $r$-subgroup for some prime $r\in\pi_a\cup\pi_b$ (and so dividing the order of $S$). Therefore, $\pi(C)$ is the union of four odd components including $\pi_a$ and $\pi_b$, and $\pi(G/G_1)$ is the fifth such component. Denote by $H$ the preimage in $G$ of a Hall $(\pi_a\cup\pi_b)$-subgroup of $S$. The group $H$ is solvable and none of its Sylow subgroups can be cyclic, which arrives us to a final contradiction.
\end{proof}

\begin{lemma}\label{directfactorRq}
Any $L_i$ isomorphic to $L$ is a direct factor in $G_1/O_3(K)$.
\end{lemma}

\begin{proof} For brevity, set $L=L_i$. Recall also that the Schur multiplier of $L$ is trivial, see Lemma~\ref{l:ree}(4). By Lemma~\ref{l:Walter}, $K/O_{2'}(K)$ is a direct factor of $G_1/O_{2'}(K)$, so we may suppose that $K$ is of odd order. Consider the upper $r$-series~\eqref{eq:r} of $K$ for $r=3$. Suppose first that $K_t=R_t=K$. Then $N=K/R_{t-1}$ is a $3'$-group. If $N_m=N$ and $N/N_{m-1}$ is the upper nontrivial factor of a chief $L$-invariant series of $N$, then by induction on the length $m$ of the series, we may assume that $L$ acts nontrivially on $N/N_{m-1}$ which is an elementary abelian $p$-group for some $p>3$. If $p^\beta$ is the $p$-exponent of $L$ (and so of $L\times L$), then Lemma~\ref{l:ReeRep} implies that $G$ contains an element of order $p^{\beta+1}$, a contradiction.

Thus, $R=R_t/K_t$ must be nontrivial. If $L\subseteq C_{(G/K_t)}(R)/R$, then $L$ is a direct factor in $G/K_t$ and we come to the situation from the previous paragraph. If $K_t=1$, then there is nothing to prove, so by induction on the order of $G$, we may assume that $V=K_t$ is an elementary abelian $p$-group for some prime $p>3$. Note that $C_K(V)\leq V$ by the definition of the upper $3$-series~\eqref{eq:r}. Let $x$ be an element in $G/K_t$ which order $p^\beta$ is equal to the $p$-exponent of $L$. Then the preimage in $G$ of the group $R\langle x\rangle$ satisfies the hypothesis of Lemma~\ref{l:hh} and provides an element of order $p^{\beta+1}$ in~$G$, a contradiction.
\end{proof}

Observe that the nonempty sets $\rho_i$, $i=3,4,5,$ and $6$, of primitive prime divisors $r_i$ of the numbers $3^{\alpha}-1$, $3^{2\alpha}-1$, $3^{\alpha}-3^{(\alpha+1)/2}+1$, and $3^{\alpha}+3^{(\alpha+1)/2}+1$, respectively (cf., Lemmas~\ref{l:bangz} and~\ref{l:GP}) are subsets of the sets $\pi_i$, $i=3,4,5,$ and $6$, respectively. Put $\rho=\rho_3\cup \rho_4\cup \rho_5\cup \rho_6$.

\begin{lemma}\label{l:alpha1} Let $L_1$ be a group over the field of order $u=3^{\alpha_1}$. Then the following hold:
\begin{itemize}
\item[{(1)}] $\alpha_1$ divides $\alpha;$
\item[{(2)}] if $r\in\rho$, then $r$ is greater than $\alpha$ and does not divide $|\Aut L_1/L_1|;$
\item[{(3)}] if $L_1\neq L$, then $(\rho_5\cup\rho_6)\cap\pi(L_1)=\varnothing$, moreover, if $\alpha_1<\alpha$, then $\rho_3$ intersects $\pi(L_1)$ trivially.    \end{itemize}
\end{lemma}

\begin{proof} (1) Note that $3^{\alpha_1}-1$ divides the exponent of $G$ which, in turn, divides $9(3^{6\alpha}-1)$. Lemma~\ref{PosLi} implies that $\alpha_1$ is odd. It follows that either $\alpha_1$ divides $\alpha$ or $\alpha_1=3\alpha$. The latter equality derives a contradiction. Indeed, if $L_1$ is a Ree group, then $(3^{\alpha_1}-3^{(\alpha_1+1)/2}+1)(3^{\alpha_1}+3^{(\alpha_1+1)/2}+1)$ does not divide $3^{6\alpha}-1$ for $L$, while if $L_1=L_2(3^{3\alpha})$, then $3^{3\alpha}-1$ does not divide the exponent of $L\times L$. Thus, $\alpha_1$ divides $\alpha$.

(2) Each $r\in\rho$ is greater than~$\alpha$ simply in view of Fermat's little theorem. None of them divides $|\Aut L_1/L_1|$, because $\alpha_1\leq \alpha$ by Item (1).

(3) Let first  $\alpha_1<\alpha$. It is clear that $(\rho_5\cup\rho_6)\cap\pi(L_1)=\varnothing$. Moreover, if $L_1=R(3^{\alpha_1})$, then $\rho_3\cap\pi(L_1)=\varnothing$, and if $L_1=L_2(3^{\alpha_1})$, then $(\rho_3\cup\rho_4)\cap\pi(L_1)=\varnothing$. If $L_1\neq L$ and $\alpha_1=\alpha$, then $L_1=L_2(3^\alpha)$ and  $(\rho_5\cup\rho_6)\cap\pi(L_1)=\varnothing$, as required.
\end{proof}

\begin{lemma}\label{kneq1}
$k=2$.
\end{lemma}

\begin{proof} If $k\neq 2$, then by above, $S=L_1\simeq R(u)$, $u=3^{\alpha_1}$, $\alpha_1$ is odd and $3\leq \alpha_1\leq \alpha$.

Suppose that $\alpha_1\neq \alpha$. Then, by Lemma~\ref{l:alpha1}(3), three primes $r_i\in\rho_i$, $i=3,5,6,$ do not divide $|L_1|$, and, by Lemma~\ref{l:alpha1}(2), do not divide $|\ov{G}|$, which contradicts Lemma~\ref{l:3dif}. Thus, $\alpha_1=\alpha$, $u=q$, and $S=L=R(q)$.

Lemma~\ref{directfactorRq} yields $\widetilde{G}_1=G_1/O_3(K)=D\times L$. Since $\rho$ intersects $\pi(\Aut{L}/L)$ trivially and $r_i\in\rho_i$, $i=3,4,5,6$ are pairwise nonadjacent in $\Gamma(L)$, at least three of them should divide the order of the solvable group~$D$. Lemma~\ref{l:criterion}(1) implies that at least two of them, say, $r$ and $s$ are adjacent in $\Gamma(D)$. If $p$ is a prime from the remaining fourth component, then $prs\in\omega(D\times L)\subseteq\omega(G)$. This contradiction completes the proof.
\end{proof}

Since $k=2$, Lemma~\ref{l:Walter} yields $\ov{G}=A_1\times A_2$, where $L_i\leq A_i\leq \Aut L_i$, $i=1,2$.

\begin{lemma}\label{l:notcyclic}
If an odd prime $r$ divides $|K|$, then a Sylow $r$-subgroup of $K$ is not cyclic.
\end{lemma}

\begin{proof}
Suppose the contrary. By Lemma~\ref{l:3dif},  there are distinct $a,b\in\{3,4,5,6\}$ such that $\pi_a\cup\pi(L_1)\neq \varnothing \neq \pi_b\cap\pi(L_2)$. Let $s\in\pi_a\cap\pi(L_1)$ and $t\in\pi_b\cap\pi(L_2)$.  We have $r\in\pi_a\cup\pi_b$, otherwise Lemma~\ref{l:cyclic2} yields $rst\in\omega(G)$, which is impossible. In particular, $r\neq 3$. The same lemma implies that if $r\in\pi_a$, then $3rt\in\omega(G)$, while if $r\in\pi_b$, then $3rs\in\omega(G)$, a final contradiction.
\end{proof}

\begin{lemma}\label{l:three2}
The set $\pi(K)$ cannot contain three primes from the different $\pi_i$, $i\in\{3,4,5,6\}$.\smallskip
\end{lemma}

\begin{proof}
If such three primes exist, then the corresponding three Sylow subgroups of $K$ are not cyclic by Lemma~\ref{l:notcyclic}. If $T$ is a Sylow $3$-subgroup of $G$, then in view of Lemma~\ref{l:cyclic1}, the solvable group $KT$ satisfies the hypothesis of Lemma~\ref{l:criterion}(2), a contradiction.
\end{proof}

\begin{lemma}\label{k=2andL=L2}
Suppose that $L_1=L_2(u)$. Then either $u$ is a power of $3$ or $u\in \{8, 13, 37\}$.
\end{lemma}

\begin{proof} Suppose the contrary. Since $u$ is not a power of $3$, the prime $3$ divides the order of a torus $T$ in $L_1$, which is equal to $(u-\varepsilon1)/(2,u-1)$ for $\varepsilon\in\{+,-\}$. Since, additionally, $u\notin \{8, 13, 37\}$, there is a prime $p>3$ also dividing the order of~$T$. Clearly, $p\in\pi_a$, where $a\in\{3,4,5,6\}$. If there is a prime $s$ from $\pi_i\cap\pi(A_2)$, where $i\in\{3,4,5,6\}\setminus\{a\}$, then $3ps\in\omega(G)$, which is impossible. It follows that $\varnothing\neq \pi(A_2)\setminus\{2, 3\}\subseteq\pi_a$.

By Lemma~\ref{l:three2}, there is $r\in\rho\setminus\pi_a$ such that $r\not\in \pi(K)$. If $r$ divides $|A_1/L_1|$, then $3r\in\omega(A_1)$, so $3rs\in\omega(G)$ for every prime $s\in\pi_a\cap \pi(A_2)$, a contradiction. Therefore, $r$ divides $|L_1|$ and does not divide $|K|\cdot|A_2|\cdot |A_1/L_1|$.

Suppose first that $r$ is not the characteristic of~$L_1$. Then a Sylow $r$-subgroup $R$ of $G$, being isomorphic to a Sylow subgroup of $L_1$, is cyclic. If $x\in R$ is an element of order~$r$, then $\pi(C_G(x))=\pi(G)$. Lemmas~\ref{l:lin}(3) and~\ref{l:LinAct} imply that the same equality $\pi(C_G(x))=\pi(G)$ holds true, if $x$ maps to a unipotent element of $L_1$, which lies in the centralizer of a field automorphism of $L_1$ that generates $A_1/L_1$. Denote by $H$ the preimage of the direct product of $L_1$ and some Sylow $3$-subgroup of $L_2$ in $G$. The group $C_G(x)\cap H$ is solvable, because $C_{L_1}(x)$ is solvable. On the other hand, $\pi(C_G(x)\cap H)$ contains $3$ and nontrivially meets $\pi_i$ for each $3\leq i\neq a$. Hence we get a contradiction in view of Lemma~\ref{l:criterion}(1).
\end{proof}

\begin{lemma}\label{OneOfIsR}
If $q>27$, then one of the $L_i$, $i=1,2$, is isomorphic to~$L$.
\end{lemma}

\begin{proof} Since $q>27$, any prime from $\rho$ does not divide the order of $L_2(u)$ for $u\in\{8,13,37\}$. So both $L_1$ and $L_2$ are groups in characteristic $3$. Suppose to the contrary that neither of $L_1$ and $L_2$  is isomorphic to~$L$. Lemma~\ref{l:alpha1} implies that $(\rho_5\cap\rho_6)\cap\pi(\ov{G})=\varnothing$. Assume that for at least one of the $L_i$, say, for $L_1$, the order of the underlying field is less than~$q$. Then, by the same Lemma~\ref{l:alpha1}, $\rho_3\cap\pi(L_1)=\varnothing$. It follows from Lemmas~\ref{l:alpha1}(2) and~\ref{l:three2} that $r\in\rho_3$ divides only the order of $L_2$. Let $H$ be the preimage in $G$ of the direct product of a Sylow $3$-subgroup of $L_1$ and a Sylow $r$-subgroup in $L_2$. Then $H$ is solvable and, as one can see, every two distinct primes from the set $\sigma=\{3, r, r_5, r_6\}$ must be adjacent in~$\Gamma(H)$. Thus, $H$ and $\sigma$ satisfy the hypothesis of Lemma~\ref{l:criterion}(2), a contradiction.

Therefore, $L_1\simeq L_2\simeq L_2(q)$. By Lemma~\ref{l:three2}, $\pi(K)\cap(\pi_3\cup\pi_4)=\varnothing$. If $3$ divides $|K|$, then a Sylow $3$-subgroup of $K$ cannot be cyclic by Lemma~\ref{l:cyclic2}. Therefore, if $R$ is a Sylow $r$-subgroup of $G$, where $r\in\pi_3\cup\pi_4$, then the solvable group $H=KR$ satisfies the conditions of Lemma~\ref{l:criterion}, a contradiction. Thus, $\pi(K)\subseteq \pi_5\cup \pi_6\cup\{2\}$. Arguing by induction on $|K|$ that $\omega(G)\not\subseteq\omega(L\times L)$, we may suppose now that $K$ is an elementary abelian $p$-group for some $p\in\pi_5\cup\pi_6$. If $x$ is an element of order $r\in\rho_3\subseteq\pi(q-1)$ from $L_1$, then $V=C_K(x)\neq 1$. This is clear if $L_1\leq C_{\ov{G}}(K)$, and follows from Lemmas~\ref{l:frob} and~\ref{l:lin}(5) otherwise. A noncyclic Sylow $3$-subgroup of $L_2$, being permutable with $x$, acts on $V$, so $3rp\in\omega(G)\setminus\omega(L\times L)$, and this contradiction completes the proof.
\end{proof}

\begin{lemma}\label{l:27}
If $L=R(27)$, then one of the $L_i$, $i=1,2$, is isomorphic to~$L$.
\end{lemma}

\begin{proof} Since $L=R(27)$, it follows that $\rho_3=\pi_3=\{13\}$, $\rho_4=\pi_4=\{7\}$, $\rho_5=\pi_5=\{19\}$, and $\rho_6=\pi_6=\{37\}$. Then $L_i\in\{L_2(8),L_2(13),L_2(27), L_2(37)\}$ for $i=1,2$.

If one of $L_i$, say $L_1$, is isomorphic to $L_2(37)$, then we arrive to a contradiction immediately, because in this case, $18\in\omega(L_1)$, $7\in\omega(L_2)$, but $18\cdot7\not\in\omega(L\times L)$. Therefore, $(\rho_5\cap\rho_6)\cap\pi(\ov{G})=\varnothing$. Further, we will mimic the proof the of previous lemma with some necessary changes.

Suppose that at least one of the $L_i$, say $L_1$, is isomorphic to $L_2(8)$. Then $\rho_3\cap\pi(L_1)=\varnothing$. Hence $r=13$ divides only the order of $L_2$ and we arrive at a contradiction, as in the first paragraph of the previous lemma. Therefore, $L_i\in\{L_2(13),L_2(27)\}$ for $i=1,2$.

If $G/G_1\neq 1$, then one of the $A_i$, say $A_1$, is a split extension of $L_2(27)$ by a field automorphism of order $3$. Then $3\cdot7\in\omega(A_1)$ and, consequently, $3\cdot 7\cdot 13\in\omega(G)$, which is impossible. Thus, $\ov{G}=S=L_1\times L_2$, in particular, the $3$-exponent of $\ov{G}$ is equal to $3$. Since the $3$-exponent of $G$ must be $9$, it follows that $3\in\pi(K)$ and we derive a contradiction as in the second paragraph of Lemma~\ref{OneOfIsR}.
\end{proof}

\begin{lemma}\label{SeqLxL}
$S=L\times L$.
\end{lemma}

\begin{proof} In view of Lemmas~\ref{OneOfIsR} and~\ref{l:27}, we may suppose that $L_1=L$. Lemma~\ref{directfactorRq} implies that $\widetilde{G}_1=G_1/O_3(K)=L_1\times D$ and $L_2$ is the unique nonabelian composition factor of~$D$. Lemma~\ref{l:alpha1}(2) yields $\rho\cap\pi(G/G_1)=\varnothing$. Assume on the contrary that $L_2\neq L$. Then Lemma~\ref{l:alpha1}(3) implies that $\pi(L_2)\cap(\rho_5\cup\rho_6)=\varnothing$. Therefore, the order of the solvable radical of $D$ is a multiple of a prime $r\in\rho_a$ for $a\in\{5,6\}$. Since Sylow $2$-subgroups of $L_2$ and $D$ itself are elementary abelian, $2r\in\omega(D)$. Let $b\in\{5,6\}\setminus\{a\}$. It follows that $2rs\in\omega(L_1\times D)\subseteq\omega(G)$ for some $s\in\rho_b$. However, $2rs\not\in\omega(L\times L)$, and we are done.
\end{proof}

\begin{lemma}\label{final}
$G=L\times L$.
\end{lemma}

\begin{proof} We claim that $G=G_1$. Indeed, if $s\in\pi(A_1/L_1)$, then there exists $a\in\{3,4,5,6\}$ such that $s\not\in\pi_a$ but $sp\in\omega(A_1)$ for some $p\in\pi_a$. Then $spr\in\omega(G)\setminus\omega(L\times L)$ for some $r\in\pi_b$, where $b\in\{3,4,5,6\}\setminus\{a\}$ and $s\not\in\pi_b$.

If $K\neq 1$, then Lemma~\ref{OneOfIsR} implies that $G/O_3(K)=L\times L\times D$ for a solvable group $D$. If  $p\in\pi(D)\setminus(\pi_5\cup\pi_6)$, then $prs\in\omega(G)$ for $r\in\rho_5$, $s\in\rho_6$. If $p\in\pi(D)\cap(\pi_5\cup\pi_6)$, then $prs\in\omega(G)$ for $r\in\pi_3$, $s\in\pi_4$.  Therefore, $D=1$ and $K$ is a $3$-group.

Let $V=K/\Phi(K)$. If $x\in L_1$ of order $r\in\rho_5$, then $U=C_V(x)>1$. This is clear if $L_1\leq C_{\ov{G}}(V)$ and follows from Lemma~\ref{l:ReeRep} otherwise. The group $L_2$ acts on $U$, so $C_V(y)>1$ for an element $y\in L_2$ of order $s\in\rho_6$ by the similar reasons. Thus, $3rs\in\omega(G)\setminus\omega(L)$, a final contradiction.
\end{proof}

\section{Proof of Theorem \ref{th2}}

Let $L=J_1$ and $G$ a finite group with $\omega(G)=\omega(L\times L)$. In view of Lemma~\ref{spectraJ1}(2), $$\mu(L)=\{6, \ 7, \ 10, \  11, \ 15, \  19\}.$$
In particular, $t(L)=4$. Thus, $G$ is nonsolvable due to Lemma~\ref{l:criterion}(2).

Let $K$ be the solvable radical of $G$ and $\overline{G}=G/K$.  By Lemma~\ref{l:Walter}, the socle $S=\Soc(\ov G)=L_1\times L_2\times \cdots \times L_k$ of $\overline{G}$ is a direct product of finite nonabelian simple groups, and every factor $L_i$, $i=1, \ldots, k$, is isomorphic to some group from the conclusion of Lemma~\ref{lem1}.

\begin{lemma}\label{l:simpleJ1}
$L_i\in\{L_2(5), L_2(11), J_1\}$ for $i=1, \ldots, k$.
\end{lemma}

\begin{proof} It follows readily from the fact that $9\not\in\mu(G)$ and $\pi(G)\subseteq\{2, 3, 5, 7, 11, 19\}$.
\end{proof}

\begin{lemma}\label{kleq2J1} $k\leq 2$.
\end{lemma}

\begin{proof} Suppose on the contrary that $k\geq 3$. Assume first that two of the factors, say $L_1$ and $L_2$, are not isomorphic to $L_2(5)$, then taking elements of orders $11$, $6$, and $5$ from $L_1$, $L_2$, and $L_3$, respectively, we obtain an element of order $330$ in $G$, which is impossible.

Suppose now that only $L_1$ is not isomorphic to $L_2(5)$. Then $k\leq 3$ and $\ov{G}=S$, because of the inequality $|\Aut L_i/L_i|\leq 2$ for $1\leq i\leq k$ and Lemma~\ref{l:Walter}.

Let $L_1\simeq J_1$. At least one of the primes $7$ and $19$ divides the order of $K$; we denote it by~$r$. Consider the upper $r$-series of $K$, see \eqref{eq:r}. We claim that $r^2\in\omega(G)\setminus\omega(L\times L)$. Indeed, if $R_t/K_t\neq 1$, then this holds true by Lemma~\ref{l:J1Rep}. Therefore, $K=K_t$, so $U=K_t/R_{t-1}\neq 1$ and $V=R_{t-1}/K_{t-1}\neq 1$. Arguing by induction on the order of $G$, we may suppose that $V$ is an elementary abelian $r$-group. Note that $C_K(V)\leq V$ by the definition of the upper $r$-series. If $L_1\subseteq C_{(G/R_{t-1})}(U)/U$, then $L_1$ is a direct factor of $G/R_{t-1}$ in view of Lemma~\ref{spectraJ1}(3), so we may apply Lemma~\ref{l:J1Rep} again. Let $x$ be an element of order $r$ in $G/R_{t-1}$. Then the preimage in $G$ of the subgroup $U\langle x\rangle$ satisfies the hypothesis of Lemma~\ref{l:hh} and provides an element of order $r^2$, as required.

Suppose that $L_1\simeq L_2(11)$. If $p\in\{7,19\}$, then $p$ divides only $|K|$ and a Sylow $p$-subgroup of $K$ is not cyclic in view of Lemma~\ref{l:cyclic2}. Assume that $11\in\pi(K)$. A Sylow $11$-subgroup of $K$ also cannot be cyclic. It follows that a Sylow $r$-subgroup of the product $H$ of a Sylow $2$-subgroup of $G$ and $K$ is noncyclic for $r\in\{2,7,11,19\}$. Now we arrive to a contradiction applying Lemmas~\ref{l:cyclic1} and~\ref{l:criterion}(2). Thus, $11\not\in\pi(K)$, so a Sylow $11$-subgroup $R$ of $G$, being isomorphic to a Sylow subgroup of $L_1$, is of order $11$. Denote by $H$ the preimage of the direct product of $L_1$ and some Sylow $2$-subgroup of $L_2$ in $G$. The group $C_G(R)\cap H$ is solvable, because $C_{L_1}(R)$ is solvable. On the other hand, $\pi(C_G(R)\cap H)$ contains $2,7,19$. Hence we get a contradiction by Lemma~\ref{l:criterion}(1).

Thus, $L_i\simeq L_2(5)$ for each $i$, so $7,11,19$ divides only $|K|\cdot|\ov{G}/S|$. As in the proof of Lemma~\ref{kLessThan3}, we derive a contradiction with Lemma~\ref{l:criterion} considering the preimage $H$ in $G$ of the normalizer in $\ov{G}$ of a Sylow $3$-subgroup of $S$ and taking its Hall $\{3, 7, 11, 19\}$-subgroup.
\end{proof}

\begin{lemma}\label{ovGeqS}
$\ov{G}=S$.
\end{lemma}

\begin{proof}
Since $|\Aut L_i/L_i|\leq 2$, it follows from Lemmas~\ref{l:Walter} and~\ref{kleq2J1}.
\end{proof}

\begin{lemma}\label{directJ1}
Any $L_i$ isomorphic to $L$ is a direct factor in $G/O_{11}(K)$.
\end{lemma}

\begin{proof} The proof is exactly the same as the proof of Lemma~\ref{directfactorRq} after replacing $r=3$ by $r=11$ and applying Lemma~\ref{l:J1Rep} instead of Lemma~\ref{l:ReeRep}.
\end{proof}

\begin{lemma}\label{LeqJ1}
One of the $L_i$ is isomorphic to $L$.
\end{lemma}

\begin{proof}
Let $\sigma=\{5, 7, 11, 19\}$. Then $L_2(5)$ and $L_2(11)$ are $D_\sigma$-groups, because $\sigma\cap\pi(L_2(5))=\{5\}$ and a Borel subgroup of $L_2(11)$ of order $55$ is its Hall $\sigma$-subgroup. Therefore, $G$ being a $D_\sigma$-group includes a solvable Hall $\{5, 7, 11, 19\}$-subgroup. An application of Lemma~\ref{l:criterion}(2) completes the proof.
\end{proof}

\begin{lemma}\label{keq2J1}
$k=2$.
\end{lemma}

\begin{proof} If $k\neq 2$, then by above $S=L=J_1$. Applying Lemma~\ref{directJ1}, we conclude that $G/O_{11}(K)=L\times D$, where $D$ is some solvable subgroup of $G/O_{11}(K)$.

Suppose first that some $p\in\{7,19\}$ does not divide~$|D|$. Then a Sylow $p$-subgroup $P$ of $G$ is of order $p$, and the centralizer $C=C_G(P)$ is solvable. On the other hand, $\pi(C)$ must contain all other prime divisors of~$|G|$. Now we get a contradiction applying Lemma~\ref{l:criterion}(1) to~$C$. Thus, $7, 19\in\pi(D)$.

Since $30\in\omega(G)$, there is $r\in\{2,3,5\}\cap\pi(D)$. In view of Lemma~\ref{l:criterion}(1), there is an element of order $ps\in\omega(D)$, where $p\neq s$ and $p, s\in\{r, 7, 19\}$. Then $11ps\in\omega(G)$, a contradiction.
\end{proof}

\begin{lemma}\label{SeqJ1xJ1}
$S=L\times L$.
\end{lemma}

\begin{proof} By Lemmas~\ref{directJ1} and~\ref{LeqJ1}, we may suppose that $G/O_{11}(K)=L_1\times D$, where $L_1=L$ and $L_2$ is the unique nonabelian composition factor of~$D$. Since $7$ and $19$ are nonadjacent in $\Gamma(L)$, at least one of them, say $r$, divides the order of $D$.  If $L_2\neq L$, then $r$ divides the order of the solvable radical $R$ of $D$. Since Sylow $2$-subgroups of $L_2$ and $D$ are elementary abelian, $2r\in\omega(D)$, so $2\cdot 7\cdot 19\in\omega(G)$, a contradiction.
\end{proof}

\begin{lemma}\label{finalJ1}
$G=L\times L$.
\end{lemma}

\begin{proof} By Lemmas~\ref{directJ1} and~\ref{SeqJ1xJ1}, we may suppose that $G/O_{11}(K)=L\times L\times D$, where $D$ is a solvable group. Suppose first that $D\neq 1$ and $p\in\pi(D)$. Take two distinct primes $r,s\in\{7, 11, 19\}\setminus\{p\}$. Then $prs\in\omega(G)\setminus\omega(L\times L)$, a contradiction.

Thus, $\ov{G}=S=L_1\times L_2$, where $L_1\simeq L\simeq L_2$ and $K$ is an $11$-group. If $K\neq 1$, we may suppose that $K$ is elementary abelian. If $L_1$ and $L_2$ acts on $K$ trivially, then $7\cdot 11\cdot 19\in\omega(G)$, which is impossible.
Suppose that one of $L_i$, say $L_1$, acts on $K$ nontrivially. Then, in view of~\cite[Theorem~1.1]{99Zal} and the fact that $11^2\not\in\omega(G)$, any chief $L_1$-invariant factor $W$ of $K$ with nontrivial action of $L_1$ is a $7$-dimensional $L_1$-module. Expecting \cite{ABCh}, one can see that given an element $x$ of order $7$ in $L_1$, the centralizer $C_W(x)$ is $1$-dimensional. Now an element $y$ of order $19$ in $L_2$ acts on $C_W(x)$ trivially. Then again $7\cdot 11\cdot 19\in\omega(G)$, a final contradiction.
\end{proof}

\medskip

\textbf{Acknowledgment.} The authors are very grateful to Maria Grechkoseeva for useful comments and fruitful discussions.

\end{document}